\documentclass[3p,twocolumn]{elsarticle}

\usepackage{geometry}
\usepackage{graphicx}
\usepackage{amsmath} 
\usepackage{amssymb}  
\usepackage{amsfonts}

\usepackage{lineno}

\usepackage{xcolor}

\usepackage{graphicx}      
\usepackage{lipsum}
\usepackage{amsmath}
\usepackage{amsfonts}
\usepackage{amssymb}
\usepackage{xcolor}
\usepackage{hyperref}
\usepackage{caption}
\usepackage{subcaption}

\newtheorem{assumption}{Assumption}
\newtheorem{theorem}{Theorem}
\newtheorem{corollary}{Corollary}
\newtheorem{proposition}{Proposition}

\newtheorem{remark}{Remark}

\journal{ISA transactions}

\begin{document}

\begin{frontmatter}


\title{Assessing the Finite-Time Stability of Nonlinear Systems by means of Physics-Informed Neural Networks}

\author[unitus]{A.~Mele}\ead{adriano.mele@unitus.it}   
\author[unina]{A.~Pironti}\ead{pironti@unina.it}

\address[unitus]{Dipartimento di Economia, Ingegneria, Società e Impresa,\\Università degli Studi della Tuscia, Viterbo, Italy}
\address[unina]{Dipartimento di Ingegneria Elettrica e delle Tecnologie dell'Informazione,\\Universit\`a degli Studi di Napoli Federico II, via Claudio 21, 80125, Napoli, Italy}

\begin{abstract}

{
In this paper, the problem of assessing the Finite-Time Stability (FTS) property for general nonlinear systems is considered. 
First, some necessary and sufficient conditions that guarantee the FTS of general nonlinear systems are provided; such conditions are expressed in terms of the existence of a suitable Lyapunov-like function. 
Connections of the main theoretical result of given in this article with the typical conditions based on Linear Matrix Inequalities (LMI) that are used for Linear Time-Varying (LTV) systems are discussed. An extension to the case of discrete time systems is also provided. 
Then, we propose a method to verify the obtained conditions for a very broad class of nonlinear systems. The proposed technique leverages the capability of neural networks to serve as universal function approximators to obtain the Lyapunov-like function. The network training data are generated by enforcing the conditions defining such function in a (large) set of collocation points, as in the case of Physics-Informed Neural Networks. 
To illustrate the effectiveness of the proposed approach, some numerical examples are proposed and discussed.
The technique proposed in this paper allows to obtain the required Lyapunov-like function in closed form. This has the twofold advantage of a) providing a practical way to verify the considered FTS property for a very general class of systems, with an unprecedented flexibility in the FTS context, and b) paving the way to control applications based on Lyapunov methods in the framework of Finite-Time Stability and Control.
}
\end{abstract}


\end{frontmatter}


\section{Introduction}
\begin{figure}[h]
    \centering
    \includegraphics[width=0.5\linewidth]{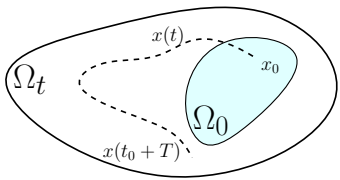}
    \caption{Graphical representation of the FTS property for a 2D system.}
    \label{fig:FTS}
\end{figure}

The concept of Finite-Time Stability (FTS), sometimes also referred to as \textit{practical stability}, was originally introduced in the Russian literature in the '50s~\cite{ftsold:1,ftsold:2,ftsold:3}; {
during the next decade, studies on this topic were carried out by the western scientific community as well~\cite{dorato:1961,weiss:1967,Michel:1972}.}
%
Roughly speaking, a dynamical system is said to be~FTS with respect to a given time-horizon~$T$, an initial time instant~$t_0$, an \emph{initial set} $\Omega_0$ and a (possibly time-varying over the time interval~$\left[t_0\,, t_0+T\right]$) \emph{trajectory set} $\Omega_t$ if, whenever the initial point of the state trajectory is contained in $\Omega_0$, the trajectory is confined inside $\Omega_t$ for all $t \in \left[t_0\,, t_0+T\right]$.\footnote{{Note that the concept of FTS used in this paper is not to be confused with another notion of finite-time stability that exists in the literature, namely the fact that a stable equilibrium is reached in a finite time interval~\cite{bhat2000finite}.}}
{(see fig.~\ref{fig:FTS}). }

FTS is a concept linked to, but independent from, Lyapunov stability, 
as it is {intrinsically} concerned with \emph{quantitative} results, substantially shifting the focus from the steady-state to the transient behaviour of the considered system.
%
However, in many cases Lyapunov-like functions can be employed to assess the FTS of a dynamical system\footnote{For brevity, we will refer to these functions as Lyapunov functions as well, even though the conditions that they are required to satisfy are slightly different from the case of standard Lyapunov stability.}; necessary and sufficient conditions for the FTS of nonlinear systems based on Lyapunov functions can be found, {for instance,} in the pioneering work by Bernfeld and Lakshmikantham~\cite{bernfeld1980practical}. {Along the same line of research, the main aim of this work is to exploit conditions similar to the ones in~\cite{bernfeld1980practical} in order to derive a practical way to assess the considered FTS property for a very broad class of nonlinear systems.}
 
In some cases, Lyapunov-like conditions for FTS can be efficiently solved by recasting them in terms of algebraic~(LMI) and Differential~(DLMI) Linear Matrix Inequalities. This fact, together with the availability of efficient tools for solving (D)LMIs, sparked a renewed interest in the topic in recent years.
This is the case, for example, of Linear Time-Varying (LTV) systems over ellipsoidal, polytopic and piecewise-quadratic domains, for which techniques based on time-varying quadratic or piecewise-quadratic Lyapunov functions have been proposed (see the book~\cite{FTSbook} and the references therein for an overview). Moreover, similar techniques can be exploited to find finite-time stabilizing controllers for the same classes of systems. 
Several results appeared in the technical literature {that extended} the FTS framework to uncertain~\cite{amato:1998,amato:2001}, hybrid~\cite{zhao2008finite,amato2013necessary}, stochastic~\cite{luan2010finite,yan2013finite} and time-delay systems~\cite{li_lyapunov_2019}, sometimes considering the case of annular domains~\cite{amato2016new,tartaglione2021annular,markovAnnular,annularStochasticLyapunov}, where both an upper and lower bound on the state variables are assigned. An approach based on the Extremum Seeking algorithm for the Finite-Time stabilization of LTV systems with unknown control direction has been recently proposed in~\cite{mele2021finite}. 
%
Some attempts have also been made in order to extend similar techniques to nonlinear systems; for example, in~\cite{amato2009sufficient} DLMI-based sufficient conditions are given for the special class of nonlinear quadratic systems. However, {for nonlinear systems, solutions in terms of (D)LMI feasibility problems are usually not available}, and different techniques must be considered. 
In this view, we aim at proposing a general approach to assess the considered FTS property for the case of nonlinear systems, which, to the best of our knowledge, is lacking in the technical literature. The main problem tackled in this article can be divided in the two following points:
\begin{itemize}
\item find necessary and sufficient conditions for the FTS of general nonlinear, non-autonomous systems; 
\item find a practical way to assess such conditions (at least for a reasonably large class of nonlinear systems for which standard LMI-based results do not apply).
\end{itemize}

To this aim, we first restate, in a slightly {simplified} form, the necessary and sufficient conditions for the FTS of continuous-time nonlinear systems originally given in~\cite{bernfeld1980practical}. {Then, we discuss how these conditions are a generalization of the ones given in~\cite{FTSbook} for linear systems over ellipsoidal domains. Furthermore, this theoretical result is also extended to the discrete-time case.} 

{Since the proposed conditions are given in terms of the existence of a suitable Lyapunov-like function, we then turn our attention to the problem of finding an expression for such function. A possible method, discussed in this work, is to design a dedicated procedure to train a suitable Neural Network (NN) in order to serve as the desired Lyapunov function.}
%
%
%
%
%
%
The use of NNs to represent Lyapunov functions is not new in the control literature, thanks to the capability of Neural Networks (NN) to serve as universal function approximators~\cite{hornik1989multilayer}.
Some applications include autonomous NL systems~\cite{abate2020formal} and robotics~\cite{chang2019neural}; in the latter paper, a learner-falsifier training procedure is proposed to speed up the training process. A similar training architecture is also used in~\cite{richards2018lyapunov}, where NNs are used to estimate the region of attraction of a given controller in the case of discrete time polynomial systems. In both~\cite{richards2018lyapunov} and \cite{gaby2021lyapunov}, multi-output NNs are used in combination with quadratic expressions to enforce the positive-definiteness of the Lyapunov functions. 
Even though NN-based function approximations are known to suffer from the so-called curse of dimensionality, in~\cite{grune2020computing} it is discussed how systems that benefit from a small-gain property have compositional Lyapunov functions which can be estimated with a number of neurons that only grows polynomially with the state dimension.
Interestingly, the training procedures used to obtain Lyapunov functions with NNs are usually aimed at enforcing the desired conditions on both the network's output and its derivatives with respect to the input parameters. As deep NNs are used to satisfy conditions involving partial derivatives by exploiting the knowledge of the equations governing a physical process, these techniques are naturally related to the framework of Physics-Informed Neural Networks (PINN)~\cite{raissi_physics-informed_2019}.

In this work, we draw ideas from these recent techniques and combine them with the general conditions for the FTS of nonlinear, non-autonomous systems with the aim to provide a flexible procedure that is capable of testing the FTS property for a broad class of dynamical systems over domains of a generic shape. 
%
To this end, we propose an approach to train a NN so that it serves as a Lyapunov function. The procedure relies on the standard Adam algorithm~\cite{kingma2015adam} to optimize a cost function designed to enforce the required conditions in a large set of collocation points. It is implemented by using off-the-shelf algorithms and can be run on GPU-enabled computers. {To the best of our knowledge, this is the first time that a practical algorithm aimed at testing the FTS property for such a broad class of systems is proposed and implemented.}

\medskip
The {rest of the} paper is organized as follows: in Section~\ref{sec:prelim}, {an overview of the notation used throughout the article} and some preliminary results are given. 
{Section~\ref{sec:mainres} presents the main results of the paper. In particular,} in section~\ref{sec:main} the main results, i.e. the necessary and sufficient conditions for {the FTS of continuous-time, nonlinear, non-autonomous systems are presented; these conditions are extended to the discrete-time case in section~\ref{subsec:discrete}}. In section~\ref{sec:LTV}, the connection between the results of section~\ref{sec:main} and existing results for continuous-time linear systems are discussed. In section~\ref{sec:NN}, a {numerical} procedure, based on NNs, to obtain a Lyapunov function {that satisfies the properties required by the main result in section~\ref{sec:main}} is {described}. In section~\ref{sec:example} some numerical examples are provided to show the effectiveness of the proposed procedure. 
{Finally, in section~\ref{sec:conclusions}, the obtained results are discussed and some concluding remarks are given}.

%

\section{Notation and Preliminaries} \label{sec:prelim}
Given a vector $x\in\mathbb{R}^n$, {$|| x ||$} indicates {its} standard euclidean norm. $\bar{\Omega}$ denotes the closure of a set $\Omega$, $|\Omega|$ its measure and $\partial \Omega$ its boundary; $f(\Omega)$ is the image of the set $\Omega$ through the function $f(\cdot)$. $M \prec 0$ indicates that the matrix $M$ is negative-definite, $M \succ 0$ that it is positive definite, $I$ denotes the identity matrix (whose dimensions will be clear from the context).

{Throughout the article, we will} consider system{s} in the form
\begin{equation} \label{eq:nlsys}
    \dot{x}(t) = f(t,x), \qquad x(t_0) = x_0 \,,
\end{equation}
where $f(t,x): J \times \mathbb{R}^n \to \mathbb{R}^n$, and $J=[t_0,T]$.
The evolution function of~\eqref{eq:nlsys} is denoted by $\varphi(t,t_0,x_0)$, while $\phi(t,t_0)$ indicates the state-transition matrix of a LTV system, for which 
\[
\phi(t,t_0) x_0 = \varphi(t,t_0,x_0) \,.
\]
{
Throughout the paper,~we will consider Lyapunov functions $V(t,x)$ that are continuous, scalar and satisfy a local Lipschitz condition with respect to $x$.  We denote the orbital derivative of such a function {a}long the trajectories of~\eqref{eq:nlsys} as (see~\cite[Chap. 1]{Yoshizawa})
\begin{equation}\label{eq:Vdotdef}
\begin{aligned}
\dot{V}(t,x) =& 
     \overline{\lim_{h\rightarrow 0^+}} \frac{1}{h} \{ V(t+h,x(t+h)) -V(t,x(t)) \} \\
     = &\overline{\lim_{h\rightarrow 0^+}} \frac{1}{h} \{ V(t+h,x+hf(t,x)) -V(t,x) \}\,.
\end{aligned}
\end{equation}
}
{Notice that, in the above definition, the upper right-hand Dini derivative appears. Usually, for a generic function of time $v(t)$, \textbf{this Dini derivative} is denoted by $D^+v(t)$; however}, to clearly differentiate between derivatives with respect to the $t$ and $x$ variables, we choose to adopt the more common dot notation $\dot{v}(t)$. When differentiable functions are considered (as in the numerical procedure described in sec.~\ref{sec:NN}), the standard and upper Dini derivatives coincide.

Finally, we adopt the following notation for compactness
\[
    \inf_{\substack{t \in J \\ x \in \Omega_t}} g(t,x) = 
    \inf_{t \in J} \inf_{x \in \Omega_t} g(t,x) \,,
\]
where the set $\Omega_t$ depends on the value of $t$.

\begin{assumption} \label{Ass:Lipschitz} The map $f(t,x)$ in equation~\eqref{eq:nlsys}
is continuous in $t$ and for each compact set $\mathcal A\subset\mathbb{R}^n$ it satisfies the Lipschitz condition
 \begin{equation}\label{eq:lipschitz_f}
     ||f(t,x)-f(t,y)|| \le \lambda_\mathcal{A}(t) ||x-y|| \qquad \forall x,y,\in \mathcal{A}\,,
 \end{equation}
where $\lambda_A(t)$ is a non-negative continuous function defined in $J$.
\end{assumption}
For each $(t^*,x^*)\in J\times\mathbb{R}^n$, assumption~\ref{Ass:Lipschitz} guarantees the existence and uniqueness of the evolution function $\varphi(t,t^*,x^*)$ for all $t$ in a suitable neighborhood of $t^*$ (see for example Theorem 3.1 in~\cite{Hale}). 
%

Given an \emph{initial set} $\Omega_0$, a (possibly time-varying) \emph{trajectory set} $\Omega_t$ and two scalars $t_0$ and $T$, system~\eqref{eq:nlsys} is said to be Finite-Time Stable (FTS) with respect to $(t_0, T, \Omega_0, \Omega_t)$ if and only if
%
\begin{equation} \label{eq:FTScond}
    x_0 \in \Omega_0 \implies x(t) \in \Omega_t \quad \forall t \in [t_0, t_0+T]\,,
\end{equation}
%
where $\Omega_0$ is a closed set, $\Omega_t$ is open and both are connected and bounded $\forall t \in J$. {Note that,} for the definition to be well-posed, it is required that $\Omega_t \supset \Omega_0$ at $t = t_0$. 





\section{Main results}\label{sec:mainres}
{In this section, the main results of this work are presented. First, necessary and sufficient conditions for the FTS of nonlinear, non-autonomous systems in the form~\eqref{eq:nlsys} under assumption~\ref{Ass:Lipschitz} are provided in thm.~\ref{thm1}. Then, these conditions are extended to the case of discrete-time systems in thm.~\ref{thm2}. The connections with the usual LMI-based conditions used for LTV systems (reported for ease of reference as thm.~\ref{thm3}) are discussed in sec.~\ref{sec:LTV}. Finally, a numerical procedure to find a closed-form expression for the required Lyapunov function, based on Neural Networks, is proposed in sec.~\ref{sec:NN}.
}

\subsection{Necessary and sufficient conditions for FTS}\label{sec:main}
The following theorem provides necessary and sufficient conditions for the FTS of system~\eqref{eq:nlsys}
\footnote{The proof closely follows the arguments in~\cite[Thm. 3.1]{bernfeld1980practical}.}.

\begin{theorem}\label{thm1}
System~\eqref{eq:nlsys} is FTS wrt $(t_0,T,\Omega_0,\Omega_t)$ if and only if there exists a
continuous function $V: J \times \mathbb{R}^n \rightarrow \mathbb{R}$, {locally Lipschitz in $x$,} such that the following conditions hold
    %
        \begin{equation}\label{eq:Vdot}
            \dot{V}(t,x) \le 0  \quad \forall t \in J, x \in \bar{\Omega}_t
        \end{equation}
        %
        \begin{equation}\label{eq:Vbound}
            \sup_{x\in\Omega_0} V(t_0,x) < \inf_{\substack{ t \in J \\ x \in \partial \Omega_t }} V(t,x)
        \end{equation}

\end{theorem}


\emph{Proof.}
We start by proving sufficiency. Assume that $\exists V: J \times \mathbb{R}^n \rightarrow \mathbb{R}$ that satisfies conditions~(\ref{eq:Vdot}-\ref{eq:Vbound}). Moreover, by contradiction, assume that, for some initial point $x_0$, $\exists t^* \in J$ such that $x(t^*)\in \partial \Omega_{t^*}$ and $x(t)\in \Omega_t$ $\forall t \in [t_0, t^*)$. 
By the definition of infimum
\[
V\left(t^*,x(t^*)\right) \ge \inf_{\substack{ t \in J \\ x \in \partial \Omega_t }} V(t,x) \,.
\]
On the other hand, eq.~\eqref{eq:Vdot} yields
\[
V\left(t^*,x(t^*)\right) \le V(t_0,x_0) \le \sup_{\xi \in \Omega_0}V(t_0,\xi) \,.
\]
which contradicts~\eqref{eq:Vbound}. Sufficiency is proven.

To prove necessity, assume that the system~\eqref{eq:nlsys} is FTS wrt $(t_0,T,\Omega_0,\Omega_t)$. We look for a continuous, locally Lipschitz function $V(t,x)$ that satisfies conditions (\ref{eq:Vdot}-\ref{eq:Vbound}). To this end consider the system
\begin{equation} \label{eq:nlsys_extended}
    \dot{x}(t) = F(t,x) \,,
\end{equation}
%
where $F$ is bounded in $J\times\mathbb{R}^n$, satisfies the global Lipschitz condition
 \begin{equation}
||F(t,x)-F(t,y)||\leq \lambda(t)||x-y||\quad \forall x,y\in\mathbb{R}^n\,,     
 \end{equation}
$\lambda(t)$ being a continuous non-negative scalar function, and finally
\begin{equation} \label{eq:f_extended}
F(t,x)=f(t,x)\,,\, \forall (t,x)\in J \times D\,,
\end{equation}
where $D$ is a compact set such that $D\supseteq\bar{\Omega}_t$ for all $t\in J$.
The existence of such a $F(t,x)$ is guaranteed by Assumption~\ref{Ass:Lipschitz}. Moreover the main result in~\cite{Derrick76} allows to conclude that the evolution function $\varphi_F(t_1,t_2,\hat x)$ of equation~\eqref{eq:f_extended} is defined, unique, and continuous in $J\times J\times\mathbb{R}^n$.
Consider
\begin{equation}\label{eq:proof_V}
    V(t,x) := \inf_{\xi \in \Omega_0} || \varphi_F(t_0,t,x) - \xi || 
\end{equation}
Note that $\varphi_F(t_0,t,x)$ is the initial point obtained by following the trajectory of~\eqref{eq:nlsys_extended} backwards from $(t,x)$ to $t_0$; the function $V(t,x)$ associates to $(t,x)$ the distance of the initial point of the trajectory from the initial set $\Omega_0$. The function $V$ is continuous and bounded on each compact subset of $\mathbb{R}^n$; we will now show that $V(t,x)$ is also Lipschitz in $x$.
First of all note that, by the Gronwall inequality, we have that
\begin{equation}\label{eq:gronwall}
||\varphi_F(t,s,x)-\varphi_F(t,s,y)||\leq \exp{\left(\int_J \lambda(\tau)\, d\tau\right)}||x-y||\,.
\end{equation}
Now consider
\begin{equation*}
    V(t,x) \le ||\varphi_F(t_0,t,x) - \xi || \quad \forall \xi \in \Omega_0 \,.
\end{equation*}
By applying the triangular inequality we obtain
\begin{equation*}
\begin{aligned}
    V(t,x) &\le || \varphi_F(t_0,t,x) - \varphi_F(t_0,t,y) || + ||\varphi_F(t_0,t,y) - \xi ||\quad \\
    &\forall \xi \in \Omega_0\,,
\end{aligned}
\end{equation*}
and this yields
\begin{equation*}
\begin{aligned}
    V&(t,x) \le \\
    &|| \varphi_F(t_0,t,x) - \varphi_F(t_0,t,y) || + \inf_{\xi \in \Omega_0} ||\varphi_F(t_0,t,y) - \xi || = \\
    &|| \varphi_F(t_0,t,x) - \varphi_F(t_0,t,y) || + V(t,y)
\end{aligned}
\end{equation*}
Similarly we have 
\[
    V(t,y) \le || \varphi_F(t_0,t,x) - \varphi_F(t_0,t,y) || + V(t,x)
\] and hence
\begin{equation}\label{eq:proof1}
    |V(t,x)-V(t,y)| \le || \varphi_F(t_0,t,x) - \varphi_F(t_0,t,y) ||\,.
\end{equation}
Now, considering equation~\eqref{eq:gronwall}, we obtain
\begin{equation} \label{eq:LipV}
|V(t,x)-V(t,y)| \le\exp{\left(\int_J \lambda(\tau)\, d\tau\right)}||x-y||\,.
\end{equation}
Equation~\eqref{eq:LipV} guarantees the existence of the orbital derivative of $V(t,x)$ as defined in~\eqref{eq:Vdotdef}. It is then trivial to verify that such derivative along the trajectories of system~\eqref{eq:nlsys_extended} is identically zero, as the initial point does not change. Now, considering equation~\eqref{eq:f_extended} and the fact that system~\eqref{eq:nlsys} satisfies the FTS property, we have
\begin{equation*}
    V(t,\varphi(t,t_0,x_0))=
    V(t,\varphi_F(t,t_0,x_0))\,, \forall (t,x_0)\in J\times\Omega_0\,,
\end{equation*}
hence,~\eqref{eq:Vdot} holds {(with the equality sign)}. 
To prove that $V(t,x)$ satisfies~\eqref{eq:Vbound}, observe that 
$\sup_{x\in \Omega_0}V(t_0,x) = 0$, 
since every point $x \in \Omega_0$ has zero distance from the set $\Omega_0$.
Now consider again a generic point $x^* \in \partial \Omega_{t^*}$, $t^*\in~J$, 
and assume that $\exists x_0^*$ s.t. $x^* = \varphi(t^*, t_0, x_0^*)$. Since system~\eqref{eq:nlsys} was assumed to be FTS, the initial point $x_0^*$ cannot be in the set $\Omega_0$. Hence, $V(t^*,x^*(t^*)) > 0$. Since this condition holds for all possible choices of~$x^* \in \partial \Omega_{t^*}$ and for all $t^* \in J$, and considering that $V(t,x)$ is continuous in $x$ and $\partial \Omega_t$ is closed and bounded, it follows that
\[
\inf_{\substack{ t \in J \\ x \in \partial \Omega_t }} V(t,x) > 0 = \sup_{x\in \Omega_0}V(t_0,x) \,,
\]
i.e. $V(t,x)$ satisfies~\eqref{eq:Vbound}.
This completes the proof.
$\hfill \blacksquare$

\begin{remark} \label{rem1}
When a finite-dimensional state space is considered, $\partial \Omega_t, \Omega_0$ will in general be compact sets, {and} the $\inf$ and $\sup$ in eq.~\eqref{eq:Vbound} {can} be replaced with $\min$ and $\max$ respectively.
\end{remark}

\begin{remark}
Notice that, differently from {classic} Lyapunov stability, in the FTS framework it is not explicitly required {that} $V(t,x)$ {is} positive-definite. 
\end{remark}

\begin{remark}
In~\cite[Thm. 3.1]{bernfeld1980practical}, a less restrictive version of~\eqref{eq:Vdot} appears, which requires that $\dot{V}(t,x) \le g(t,V(t,x))$ for $(t,x) \in J \times \bar{\Omega}_t$, with $g:\mathbb{R}\times\mathbb{R} \to \mathbb{R}^+$. With this choice, the orbital derivative of $V$ is allowed to take positive values, but the price for this increased freedom is that condition~\eqref{eq:Vbound} must be complicated by resorting to the comparison lemma~\cite[Lemma 3.4]{khalil}. 
However, there is no indication on how to choose the function $g(t,V(t,x))$ \emph{a priori}, and even the necessity proof in~\cite[Thm. 3.1]{bernfeld1980practical} eventually assumes $g \equiv 0$. 
For these reasons, we set $g \equiv 0$ from the beginning; this choice also simplifies the training procedure of sec.~\ref{sec:NN}. 
\\
Moreover, with respect to~\cite[Thm. 3.1]{bernfeld1980practical}, the formulation of thm.~\ref{thm1} automatically takes into account the possibility that the trajectory domain $\Omega_t$ is time-varying. This, in turn, makes condition (\textit{b-ii}) in~\cite[Thm. 3.1]{bernfeld1980practical} unnecessary.
\end{remark}


In view of the numerical implementation of the conditions in thm.~\ref{thm1}, it is useful to state the following corollary.

\begin{corollary}\label{cor1}
System~\eqref{eq:nlsys} is FTS wrt $(t_0,T,\Omega_0,\Omega_t)$ if and only if there exists a continuous function $V: J \times \mathbb{R}^n \rightarrow \mathbb{R}$ such that condition~\eqref{eq:Vdot} holds and condition~\eqref{eq:Vbound} is replaced by
%
        \begin{equation}\label{eq:Vbound-2}
         \inf_{\substack{ t \in J \\ x \in \partial \Omega_t }} V(t,x) - \sup_{x\in\Omega_0} V(t_0,x) > \alpha
        \end{equation}
for some fixed value of $\alpha>0$.
\end{corollary}

\emph{Proof.}
Suppose that~\eqref{eq:Vbound} holds with
\[
\inf_{\substack{ t \in J \\ x \in \partial \Omega_t }} V(t,x) - \sup_{x\in\Omega_0} V(t_0,x) > \varepsilon
\]
for some $V(t,x)$ and $\varepsilon > 0$. 
The function $\tilde{V}(t,x) := \beta V(t,x)$, with $\beta > 0$ satisfies condition~\eqref{eq:Vdot}, since $\dot{\tilde{V}}(t,x) = \beta \dot{V}(t,x) \le 0$. Choosing $\beta = {\alpha}/{\varepsilon} $, it is now easy to verify that $\tilde{V}(t,x)$ satisfies condition~\eqref{eq:Vbound-2} as well. 

%
$\hfill \blacksquare$

{The proof of this corollary is trivial, and it stems from the observation that a Lyapunov function can be rescaled by an arbitrary constant. However, the usefulness of this result is in the fact that it allows to freely choose the parameter $\alpha$, which becomes a tuning parameter for the numerical algorithm described in sec.~\ref{sec:NN}. Moreover, it can be used to draw connections between well-assessed results on the FTS of LTV systems and thm.~\ref{thm1}; this is done in sec.~\ref{sec:LTV}.}

\subsection{Extension to discrete-time systems}\label{subsec:discrete}

Consider a system in the form
\begin{equation}\label{eq:dtsys}
    x(k+1) = f(k,x(k)) \,.
\end{equation}
System~\eqref{eq:dtsys} is said to be FTS with respect to $(\Omega_0, \Omega_k, N)$ if and only if, by definition, $x(0) \in \Omega_0 \implies x(k) \in \Omega_k$ for $k = 0, 1, ..., N$ 
(notice that there is no loss of generality in assuming that the index $k$ starts from $0$). As for the continuous time case, we assume $\Omega_0$ closed and $\Omega_k$ open. Moreover, we define the set 
%
\[
    \Omega_k^{fwd} := [f(k,\bar{\Omega}_k) \setminus \Omega_{k+1} ] \cup \partial f(k,\bar{\Omega}_k)  
\]
%
%
and
\[
dV(k,x) := V(k+1,f(k,x))-V(k,x) \,.
\]
The following theorem provides necessary and sufficient conditions for the FTS of~\eqref{eq:dtsys}.

\begin{theorem}\label{thm2}
System~\eqref{eq:dtsys} is FTS wrt $(\Omega_0, \Omega_k, N)$ if there exists a function $V(k,x)$ continuous in $x$, $V: \{0,...,N\} \times \mathbb{R}^n \rightarrow \mathbb{R}$ such that the following conditions hold
%
    \begin{equation}\label{eq:td_i}
        dV(k,x) \le 0 \,, \quad x\in\Omega_k \quad k = 0, ..., N-1
    \end{equation}
    \begin{equation}\label{eq:td_ii}
        \sup_{x\in\Omega_0} V(0,x) < \inf_{\substack{0\le k < N \\ x \in \Omega^{fwd}_k}} V(k+1,x)\,.
    \end{equation}
%
Moreover, if $f(k,x)$ is continuous with respect to $x$ and {it admits an inverse $f^{-1}(k,x)$ for $k \in \{0,...,N\}$}, then conditions (\ref{eq:td_i}-\ref{eq:td_ii}) are also necessary.
\end{theorem}

\emph{Proof.}
To prove sufficiency, assume that $\exists V:\{0,...,N\} \times \mathbb{R}^n \to \mathbb{R}$ that satisfies conditions (\ref{eq:td_i}-\ref{eq:td_ii}). Moreover, assume that, for some initial point $x_0$, $\exists k^* \in \{0,...,N-1\}$ such that $x(k^*)\in \Omega_{k^*}$ and $x(k^*+1) \notin \Omega_{k^*+1}$.
Since $x(k^*+1) \in f(k,\Omega_k) \subset f(k,\bar{\Omega}_k)$, it follows that $x(k^*+1) \in \Omega^{fwd}_k$. By definition, we have that
\begin{equation}
    V\left(k^*+1,x(k^*+1)\right) \le \inf_{\substack{0\le k < N \\ x \in \Omega^{fwd}_k}} V(k+1,x) \,.
\end{equation}
On the other hand, condition~\eqref{eq:td_i} ensures that 
%
\[
    V\left(k^*+1,x(k^*+1)\right) \le V(0,x_0) \le \sup_{x\in\Omega_0} V(0,x) \,,
\]
%
which contradicts condition~\eqref{eq:td_ii}.

To prove necessity, arguments similar to those used in thm.~\ref{thm1} can be used. In particular, it can be shown that the candidate Lyapunov function
\[
    V(k,x) = \inf_{\xi \in \Omega_0} || \varphi(0,k,x) - \xi ||
\]
satisfies the conditions of the theorem if the considered system is FTS. However, for $V(k,x)$ to be well-defined, we need the function $f(k,x)$ to be invertible for all $k$.
%
%
$\hfill \blacksquare$


\subsection{Connection with continuous-time LTV systems over ellipsoidal domains}\label{sec:LTV}


The definition of FTS can be recast in terms of ellipsoidal initial and trajectory domains by substituting equation~\eqref{eq:FTScond} with
\begin{equation}
    x_0^T R x_0 \le 1 \implies x(t)^T \Gamma(t) x(t) < 1 \quad \forall t \in J
\end{equation}
where $R$ is a positive-definite, symmetric matrix and $\Gamma(t)$ is a positive-definite, symmetric, matrix-valued function. The well-posedness condition becomes $R \succ \Gamma(t_0)$.
For the FTS of LTV systems, several equivalent necessary and sufficient conditions are given in~\cite[Thm.~2.1]{FTSbook}. In particular, the following theorem is proven.

\begin{theorem}~\cite[Thm. 2.1]{FTSbook}\label{thm3}
    Consider the LTV system
    \begin{equation}\label{eq:ltvsys}
        \dot{x}(t) = A(t)x(t), \quad x(t_0) = x_0
    \end{equation}
    with the initial domain $\Omega_0 = \{x | x^T R x \le 1\}$ and the trajectory domain $\Omega_t := \{x | x^T \Gamma(t) x < 1\}$, where $\Gamma(t_0)\prec R$.
    System~\eqref{eq:ltvsys} is FTS wrt $(t_0,T,\Omega_0,\Omega_t)$ if and only if there exists a positive-definite matrix-valued function $P(\cdot)$ that satisfies the following DLMI conditions
    %
    %
    \begin{subequations}\label{eq:FTS-LTV}
    \begin{align}
      \dot{P}(t) + P(t)A^T(t) + A(t)P(t) \prec  0\,, t \in J \label{eq:FTS-LTV1}
    \\
      P(t)  \succ  \Gamma(t)\,, t \in J \label{eq:FTS-LTV2}
    \\
      P(t_0)  \prec R \label{eq:FTS-LTV3}
    \end{align}
    \end{subequations}
\end{theorem}

%
This theorem can be shown to be a particular case of thm.~\ref{thm1} {in sec.~\ref{sec:main}}. 
Let us start by choosing $V(t,x) = x(t)^T P(t) x(t)$ (the validity of this choice will be discussed at the end of this section). Then, if a strict inequality is considered\footnote{The equality case can be linked to the DLE condition (2.3a) in~\cite{FTSbook}. This point is not discussed here for brevity.}, property~\eqref{eq:Vdot} reduces to~\eqref{eq:FTS-LTV1}, while condition~\eqref{eq:Vbound} can be recovered by observing that~(\ref{eq:FTS-LTV2}-\ref{eq:FTS-LTV3}) yield:
\begin{subequations}
\begin{align}
    V(t,x) > x(t)^T \Gamma(t) x(t) \implies \inf_{\substack{ t \in J \\ x \in \partial \Omega_t}}{V(t,x)} > 1 \label{eq:ii'}
    \\
    V(t_0,x) < x(t_0)^T R x(t_0) \implies \sup_{x\in \Omega_0}{V(t_0,x)} < 1 \label{eq:ii''}
\end{align}
\end{subequations}
%
%
Hence, conditions~(\ref{eq:ii'}-\ref{eq:ii''}) imply~\eqref{eq:Vbound}.
On the other hand, once we fix $V(t,x) = x^T P(t) x$ the inverse implication is also true.
To show this, first of all observe that 
$V(t,x)$ can always be scaled up in such a way that
\[
    \inf_{\substack{ t \in J \\ x \in \partial \Omega_t}}{x^T P(t) x} > 1 > \sup_{x\in \Omega_0}{x^T P(t_0) x} \,,
\]
 without affecting the condition $\dot{V}\le0$ (cfr. cor.~\ref{cor1}), i.e. conditions~(\ref{eq:ii'}-\ref{eq:ii''}) hold. 
 The first inequality can be rewritten as
\[
    x^T P(t) x > x^T \Gamma(t) x \qquad \forall x: x^T \Gamma(t) x = 1 \,.
\]
But since $\Gamma\succ0$, for any other choice of $x\in\mathbb{R}^n$ there exists a scalar $c$ such that $x^T \Gamma(t) x = c^2$. This means that the point $\frac{x}{c}$ lies on $\partial \Omega_t$, and hence
\[
  \frac{1}{c^2} x^T \Gamma(t) x = 1 < \frac{1}{c^2} x^T P(t) x \implies x^T \Gamma(t) x < x^T P(t) x  
\]
i.e. $P(t) \succ \Gamma(t)$.
With a similar argument it can be shown that $P(t_0) \prec R$, since
\[
\sup_{x\in\Omega_0}V(t_0,x)<1 \implies \sup_{x\in\partial\Omega_0}V(t_0,x)<1 \,.
\]

To conclude this section, we discuss how the choice $V(t,x) = x^T(t) P(t) x(t)$ derives from the same arguments in the necessity part of the proof of thm.~\ref{thm1}. 
In {thm.}~\ref{thm1} we considered a function $V(t,x)$ as in~\eqref{eq:proof_V},
i.e. the distance between the initial point of the trajectory and {the set}~$\Omega_0$. We observed that, whenever $x(t) \in \partial \Omega_t$, such distance must be strictly positive if the system is FTS. In the case of ellipsoidal domains, this is equivalent to 
\begin{equation}\label{eq:ltv_cond1}
    x^T(t)\Gamma(t)x(t) = 1 \implies x_0^T R x_0 > 1 \,.
\end{equation}
Using the state-transition matrix $\phi(t,t_0)$ of~\eqref{eq:ltvsys}, we have
\[
    x^T(t)\Gamma(t)x(t) = x_0^T \phi(t,t_0)^T \Gamma(t) \phi(t,t_0) x_0(t) \,,
\]
and hence~\eqref{eq:ltv_cond1} reduces to
%
\[
    \phi(t,t_0)^T \Gamma(t) \phi(t,t_0) \prec R \,.
\]
%
Moreover, since $\phi(t,t_0)^{-1} = \phi(t_0,t)$
\begin{equation}\label{eq:ltv_cond3}
     Q(t) := \phi(t_0,t)^T R \phi(t_0,t) \succ \Gamma(t) \,,
\end{equation}
which closely resembles~\eqref{eq:FTS-LTV2}. Since $\phi(t_0,t_0) = I$, {we immediately find that} 
\[
Q(t_0) = R
\]
(compare these conditions with~\cite[Thm. 2.1 (ii-iii)]{FTSbook}). Finally, by continuity we can choose an arbitrarily small $\varepsilon > 0$ such that conditions~\eqref{eq:FTS-LTV2}-\eqref{eq:FTS-LTV3} hold for 
\[
P(t) := Q(t)(1-\varepsilon) \,.
\]
Observe that, at its core, this choice of $Q(t)$ again reduces to following the trajectories of~\eqref{eq:ltvsys} backwards and evaluating the distance of the initial point $(t_0,x_0)$ from the initial set through the quadratic forms associated to $Q(t)$ and $R$.

\subsection{Lyapunov function approximation through Physics-Informed Neural Networks}\label{sec:NN}

In this section, we propose a method to train a NN so that its output provides a Lyapunov function for continuous time systems which are FTS.
%
Consider the following cost functional
\begin{equation}
    L(\dot{V},V,\delta_2) = L_1(\dot{V}) + L_2(V,\delta_2) \,,
\end{equation}
where 
\begin{subequations}\label{eq:loss}
\begin{align}
    L_1 = 
    \frac{1}{T}
    \int_{J}
    \frac{1}{|\bar{\Omega}_t|}\int_{\bar{\Omega}_t} 
    \left( 
    \max\{\dot{V}(\tau,x),0\}
    \right)^2 dx d\tau \label{eq:loss1}
\\
    L_2 = 
    \frac{1}{T}
    \int_{J}
    \frac{1}{|\partial\Omega_t|}\int_{\partial\Omega_t} 
    \left( 
    \max\{\delta V_{b}(\tau,x) + \delta_2,0\}
    \right)^2 dx d\tau  \,.   \label{eq:loss2}
\end{align}
\end{subequations}
In~\eqref{eq:loss2}, the parameter $\delta_2 > 0$ is introduced to enforce the strict inequality in~\eqref{eq:Vbound} (see discussion below). 
The term $\delta V_{b}$ appearing in~\eqref{eq:loss2} is defined as
\[
    \delta V_{b}(t,x) = \sup_{\xi \in \Omega_0} \{V(t_0,\xi) \}
    - V(t,x)\,,
\]
and is evaluated on $\partial \Omega_t$. 
%
%
%
The following result shows that the problem of finding a Lyapunov function for a system in the form~\eqref{eq:nlsys} is equivalent to finding a $V(t,x)$ such that $L(\dot{V},V,\delta_2) = 0$ for some (arbitrary) value of $\delta_2$.

\begin{proposition}
    A continuous function $V(t,x)$ satisfies the conditions of thm.~\ref{thm1} if and only if $L(\dot{V},V,\delta_2) = 0$ for some value of $\delta_2>0$.
\end{proposition}

\emph{Proof.}
If condition~\eqref{eq:Vdot} holds, then $\dot{V}(t,x) \le 0$ for every $t \in J$ and for every $x \in \Omega_t$. Hence, $\max\{\dot{V}(\tau,x),0\} = 0$ and $L_1(\dot{V}) = 0$. On the other and, since $L_1(\dot{V})$ is defined as the integral of a non-negative quantity, $L_1(\dot{V})=0$ implies that its argument is equal to zero for all $t \in J$ and $x \in \Omega_t$, i.e. that $\dot{V} \le 0$. 
Similarly, $L_2(V,\delta_2) = 0$ if and only if $\delta V_b(t,x,\delta_2)\le 0$ for every $t \in J$ and for every $x \in \partial\Omega_t$. But this is equivalent to requiring that \begin{equation}
\begin{aligned}
     0 &\ge 
     \sup_{\substack{t \in J \\ x \in \partial \Omega_t}}  \delta V_b(t,x) + \delta_2 \\
     & = \sup_{\substack{t \in J \\ x \in \partial \Omega_t}}\left\{ 
     \sup_{\xi \in \Omega_0} \{V(t_0,\xi) \}
        - V(t,x) \right\} + \delta_2 \\
    & = \sup_{\xi \in \Omega_0} \{V(t_0,\xi) \} - \inf_{\substack{t \in J \\ x \in \partial \Omega_t}}
    \{ V(t,x) \} + \delta_2 \,.
\end{aligned}
\end{equation}
In view of corollary~\ref{cor1}, if this condition is satisfied for some value of $\delta_2>0$, then also condition~\eqref{eq:Vbound} holds.
The proof is completed by observing that $L~=~0$ is
equivalent to $L_1~=~0~\wedge~L_2~=~0$, since $L_{1,2} \ge 0$ by definition.
$\hfill \blacksquare$

In practice, 
we restrict our attention to functions $V(t,x)$ of class $C^1$ by choosing continuously differentiable activation functions. The universal approximation property of NNs ensures that, even if $C^1$ activation functions are chosen, any continuous function can be approximated with arbitrary precision by the network, provided that a sufficiently large number of nodes is employed.
It is worth to notice that it is usually required that the function approximated by the NN is evaluated on compact subsets of $\mathbb{R}^n$, which is precisely the case in FTS, where the conditions on $V(t,x)$ must be satisfied only on the (compact) domains of interest.
%
For $V(t,x) \in C^1$, the orbital derivative~\eqref{eq:Vdotdef} can be evaluated as
\begin{equation}\label{eq:Vdot_nn}
    \dot{V}(t,x) =
     \frac{\partial V}{\partial t} + \mathcal{L}_f V =
     \frac{\partial V}{\partial t} + \frac{\partial V}{\partial x}  f(t,x) \,,
\end{equation}
{where $\mathcal{L}_fV$ denotes the Lie derivative of $V(t,x)$ along $f(t,x)$.}
The partial derivatives of $V$ with respect to $x$ and $t$ can be computed analytically via automatic differentiation~\cite{baydin_automatic_nodate}.
To train the network, we sample the cost function $L(\dot{V},V)$ in a (large) set of collocation points. Let us define
%
%
\begin{equation}\label{eq:apploss1}
    \hat{L}_1(\dot{V},\delta_1) = \frac{1}{N_c} \sum_{i=1}^{N_c}
    \left(
    \max\{\dot{V}(x_i,t_i)+\delta_1,0\}
    \right)^2 \\,,
\end{equation}
where $N_c$ is the total number of considered collocation points in $t \in J, \, x \in \bar{\Omega}_t$, and
\begin{equation}\label{eq:apploss2}
\begin{aligned}
    & \hat{L}_2(V,\delta_2) = \frac{1}{N_b} \sum_{i=1}^{N_b}
    \left(
    \max\left\{\delta \tilde{V}_b(t_i,x_i) + \delta_2, 0\right\}
    \right)^2
    \\
    & \delta \tilde{V}_b(t_i,x_i) = \max_{j \in [1, N_0]} \{V(t_0,\xi_j)\} - V(t_i,x_i) \,,
\end{aligned}
\end{equation}
where $N_b$ and $N_0$ are the numbers of points taken in $t \in J, \, x \in \partial \Omega_t$ and in $\Omega_0$ respectively. The $N_c$ internal collocation points are obtained by random sampling in both time and space, while the $N_b$ boundary points are computed by discretizing the time interval in $N_t$ discrete instants and then taking $n_b = N_b/N_t$ points for each time instant. The $N_0$ points are uniformly distributed in the initial domain $\Omega_0$.
In eq.~\eqref{eq:apploss1} we introduced a tolerance parameter $\delta_1$ to penalize the points where $\dot{V}$ in $\hat{L}_1$ is exactly equal to or slightly smaller than zero, similarly to the $\delta_2$ parameter in eqns.~\eqref{eq:loss2}-\eqref{eq:apploss2}. Moreover, in view of remark~\ref{rem1}, we used the maximum in place of the supremum in the definition of $\hat{L}_2$.
The resulting (approximate) loss function is
\begin{equation}
    \hat{L}(\dot{V},V,\delta_1,\delta_2) = \alpha_1 \hat{L}_1(\dot{V},\delta_1) + \alpha_2 \hat{L}_2(V,\delta_2) \,,
\end{equation}
where the weights $\alpha_{1,2} > 0$ have been introduced as additional tuning parameter{s in} the algorithm.

The network is trained by using the matlab implementation of the standard Adam algorithm~\cite{kingma2015adam}. The default parameters have been used, i.e. a global learning rate of $0.001$, a gradient decay factor of $0.9$ and a squared gradient decay factor of $0.999$. {This choice proved to be enough to obtain good convergence properties in the considered examples.} 
To compute {$\hat{L}_1$}, the $N_c$ internal collocation points are divided into mini-batches, while all the $N_b + N_0$ points on $\partial \Omega_t$ and in the initial domain $\Omega_0$ are used at each training iteration to compute $\hat{L}_2$. 
This choice reflects the fact that~\eqref{eq:Vbound} only requires that the infimum of {$V(t,x)$} over all the boundary points is larger than {the} supremum {of $V(t_0,x)$} over the initial domain. {H}owever, in the definition of $\hat{L}_2$ we considered all the {collocation} points on $\partial \Omega_t$ to obtain a well-behaved estimate of the cost function gradients with respect to the network parameters.
The training procedure terminates when $\dot{V}(t,x) \le 0$ over all the $N_c$ {internal} collocation points and $\delta{\tilde{V}}_b(t,x) < 0$ over all the $N_b$ {boundary} collocation points. 
Finally, since we are resorting to an approximate version of the cost function $L(\dot{V},V,\delta_2)$, after a solution is found we verify that the theorem conditions are satisfied over a set of test points different from the collocation points used in the training procedure.

Before concluding this section, it is worth to observe that the numerical implementation of the conditions of theorem~\ref{thm1} makes them only necessary in practice. If the considered system is not FTS, then the procedure will usually not converge to a feasible solution; however, it may happen that the algorithm does not converge even if the system is FTS, for instance when the number of nodes in the NN is not large enough to correctly represent the desired Lyapunov function.

\section{Numerical examples}\label{sec:example}

\begin{figure*}[h]
    \centering
    \includegraphics[width=0.8\linewidth]{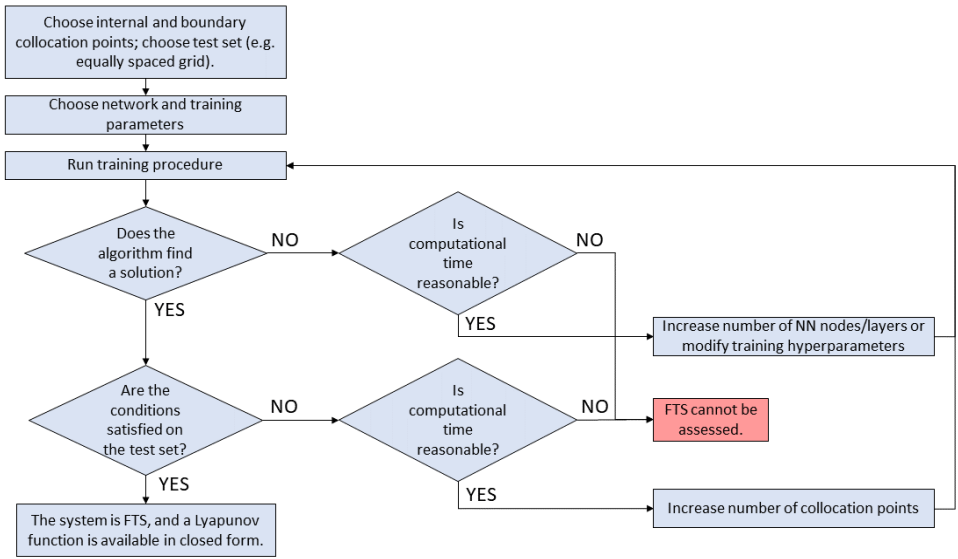}
    \caption{Flow diagram describing how the proposed technique has been applied to the examples of sec.~\ref{sec:example}.}
    \label{fig:flowchart}
\end{figure*}

In this section we illustrate the proposed approach through some numerical examples.

It is worth to observe that the conditions in thm.~\ref{thm1} are based on the orbital derivative of the function $V(t,x)$, which can be evaluated by means of automatic differentiation techniques when $V$ is expressed through a neural network. On the other hand, the conditions in thm.~\ref{thm2} are based on recursive equations, and require the computation of the set $\Omega_k^{fwd}$ at each time step. In the following, we will exploit the availability of automatic differentiation routines implemented in Matlab, and only focus on continuous-time systems. 
Note that, for simplicity, in the proposed examples ellipsoidal domains are considered. However, the algorithms only uses the collocation points placed in such domains and on their boundary, hence the technique can be easily {applied} to domains of generic shape.

{
The proposed procedure has been applied to the considered examples as shown in fig.~\ref{fig:flowchart}. For all the examples, we used a Multi-Layer Perceptron architecture with \emph{softplus} activation layers (the parameters for each example are reported in table~\ref{tab2}). After an initial choice of the collocation point and of the NN and training parameters, the optimization procedure is run. If a solution is found, the required conditions are further tested on a suitable set of test points (e.g. a finely spaced grid over the considered domain). If the test succeeds, the property assessment can be considered to be successful. If the test fails, the training set of collocation points can be enlarged, e.g. by including the points where the test fails, and the procedure can be run again. If the algorithm fails to converge, the NN and training parameters can be modified (e.g. increasing the number of nodes/layers of the NN to make it more expressive, or lowering the learning rates to improve stability) before running the procedure again. If no solution is found with a reasonable set of parameters and/or in a reasonable time, the test is considered to have failed (red box in fig.~\ref{fig:flowchart}). This could be due to the fact that the system does not satisfy the FTS property or to numerical reasons.
It is worth to notice that this characteristic is common to most Lyapunov techniques, which often do not provide a straightforward way to verify that the system is unstable (in the classic or Finite-Time sense). 
}
%

All the examples in this section have been implemented in matlab and run on a laptop equipped with an Intel Core i7-8750H CPU @2.20GHz processor, 16 GB of RAM and an NVIDIA Geforce GTX 1050 Ti graphic card. The code is available at \url{https://github.com/AdrianoMele/NeuralNetFTS}.

\subsection*{Example 1}
We start by validating the procedure on a toy example, i.e. a 2D LTI system with fixed ellipsoidal domains. 
In particular, consider the system
\begin{equation}
    \dot{x}(t) = -
    \begin{bmatrix}
        0.1 \quad 0 \\ 0 \quad 0.1
    \end{bmatrix}
    x(t) \,,
\end{equation}
with $t_0 = 0, \, T = 1$ and the domains defined by
\begin{equation}
\begin{aligned}
    &\Omega_0 = \{ x \in \mathbb{R}^2 | x^T R x \le 1 \}, \, R = 0.3I
    \\
    &\Omega_t = \{ x \in \mathbb{R}^2 | x^T \Gamma x \le 1 \}, \, \Gamma = 0.25I \,.
\end{aligned}
\end{equation}
{Classical FTS results ensure that the solution for this example can be found in terms of a time-invariant quadratic Lyapunov function $V_Q(x) = x^T(t) P x(t)$ (also see the discussion in section~\ref{sec:LTV}). In fact, the system is also Lyapunov-stable, and it can be shown that, in this case, the FTS conditions in thm.~\ref{thm3} can be reduced to the problem of finding a static Lyapunov function that also satisfies conditions~(\ref{eq:FTS-LTV2}-\ref{eq:FTS-LTV3}).
We apply the technique proposed in section~\ref{sec:NN} to the considered system and domains; we expect the algorithm to find a solution $V(x)$ that is very close to a quadratic function.}

\begin{figure}
     \centering
     \begin{subfigure}[t]{0.8\linewidth}
         \centering
         \includegraphics[width=0.7\linewidth]{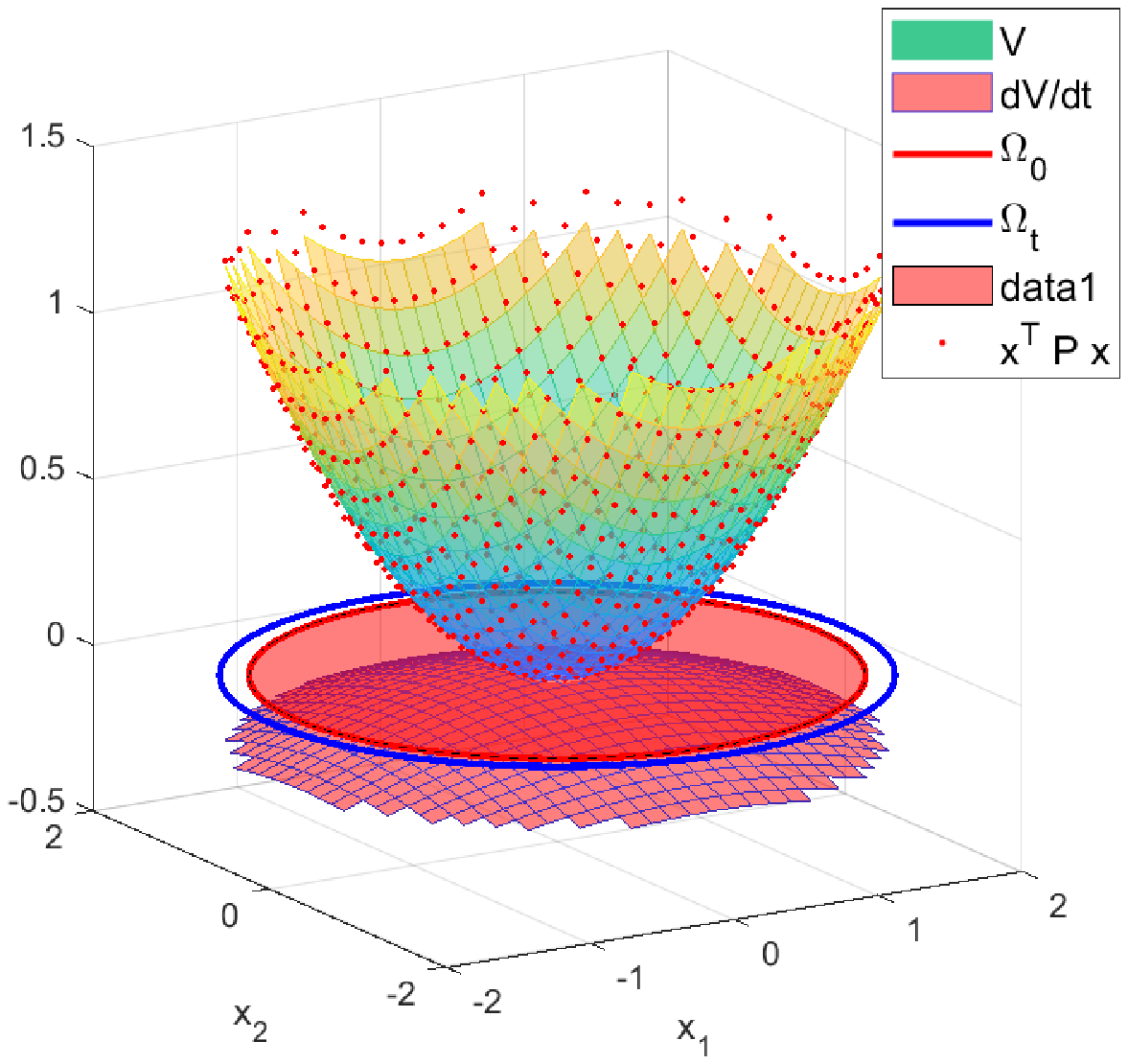}
         \caption{}
    \label{fig:ex1_1}
     \end{subfigure}
     \hfill
     \begin{subfigure}[t]{0.8\linewidth}
         \centering
         \includegraphics[width=0.7\linewidth]{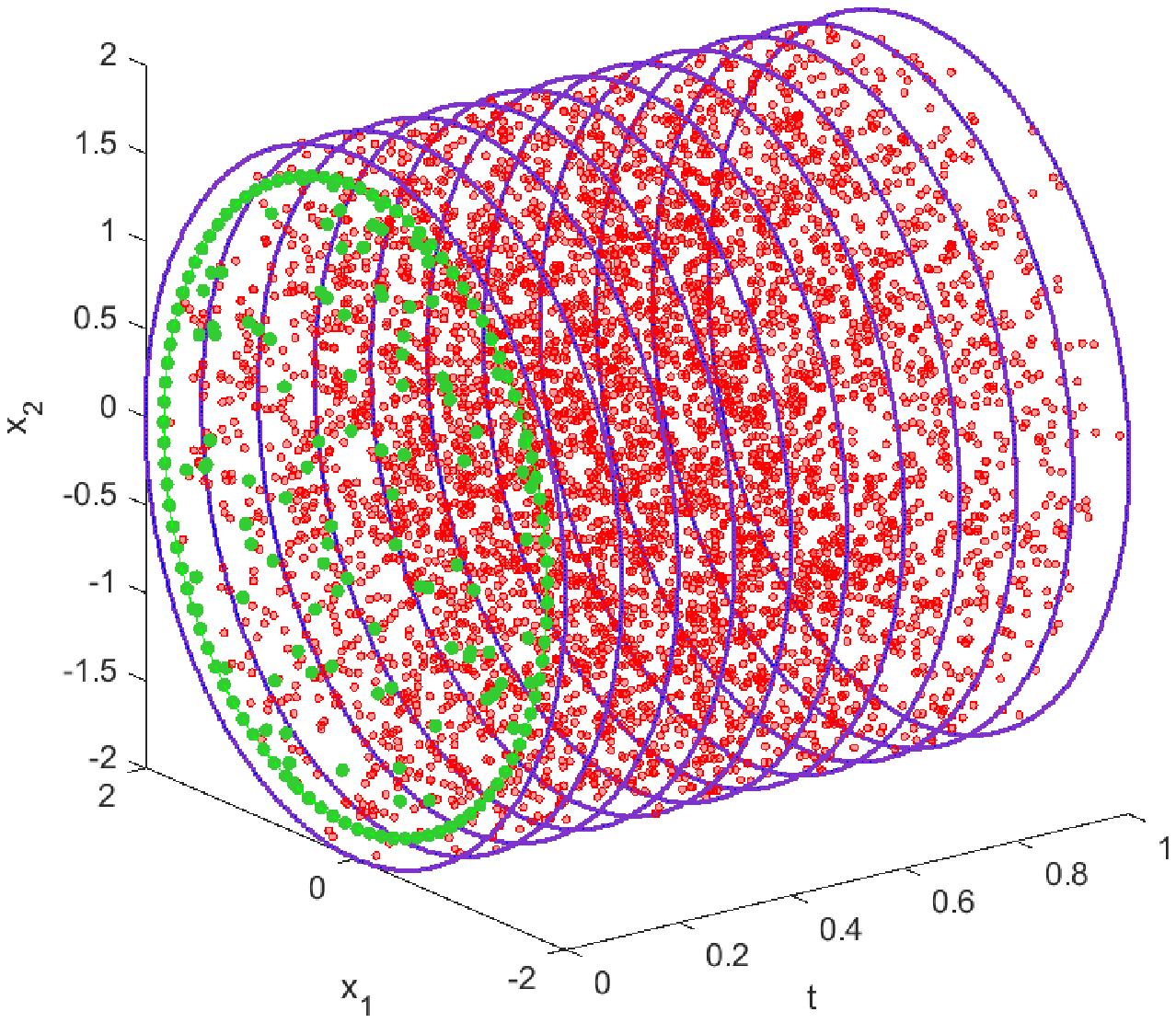}
         \caption{}
    \label{fig:ex1_2}
     \end{subfigure}
    \caption{(a) Resulting $V(t,x)$ and $\dot{V}(t,x)$ for example \#1; the red dots show the quadratic approximant (a constant has been added to $V(x)$ so that $V(0) = 0$ for graphical reasons). 
    (b) Collocation points for example \#1. The green dots show the $N_0$ points in the initial domain, the red dots are the $N_c$ internal points used to enforce~\eqref{eq:Vdot} and the blue dots are the $N_b = n_b N_t$ boundary points used to enforce~\eqref{eq:Vbound}.}
\end{figure}

The training parameters used in this example are reported in tables~\ref{tab1}-\ref{tab2}. The algorithm is able to find a solution after only 2 epochs (about $2$~s of training). The resulting $V(x)$ is shown in fig.~\ref{fig:ex1_1}, while the collocation points {and the initial and trajectory domains} are shown in fig.~\ref{fig:ex1_2}. 
A quadratic best-fit of $V(x)$ (red dots in fig.~\ref{fig:ex1_1}) results in the symmetric positive-definite matrix
%
\[
    P = 
    \begin{bmatrix}
        0.3064 \quad -0.0063 \\ -0.0063 \quad 0.3120
    \end{bmatrix}
\]
The eigenvalues of $PA^T + AP$ are $[-0.06322, \, -0.06044]$, i.e. the resulting $P$ satisfies condition~\eqref{eq:FTS-LTV1}. On the other hand, the eigenvalues of $P$ are $[0.3022, \, 0.3161]$, i.e. condition~\eqref{eq:FTS-LTV3} does not hold. As discussed in sec.~\ref{sec:LTV}, a function $V_Q(x)$ that satisfies~(\ref{eq:FTS-LTV2}-\ref{eq:FTS-LTV3}) can be obtained by a suitable scaling of $P$.




\subsection*{Example 2} 
Consider the following system~\cite[(1.1.5)]{lakshmikantham_book}
\begin{equation}\label{eq:laksh_ex}
\begin{aligned}
    \dot{x_1}(t) &= -x_1 -x_2 +k(x_1-x_2)(x_1^2+x_2^2) \\
    \dot{x_2}(t) &= x_1 -x_2 +k(x_1+x_2)(x_1^2+x_2^2) \\
\end{aligned}
\end{equation}
where the constant $k$ is chosen equal to $0.1$. {This system provides an interesting analytical benchmark for the considered numerical procedure. In fact,} in~\cite{lakshmikantham_book} it is shown that, letting $r^2(t)=x_1^2(t)+x_2^2(t)$, the solutions of~\eqref{eq:laksh_ex} have the property
\begin{equation}\label{eq:laksh_bound}
    r^2(t) = \frac{1}{k\mu}r_0^2 \,, \quad 
    \mu = r_0^2 + \left(\frac{1}{k} - r_0^2\right) e^{2(t-t_0)}
\end{equation}
 
\begin{figure}
    \centering
    \includegraphics[width=\linewidth]{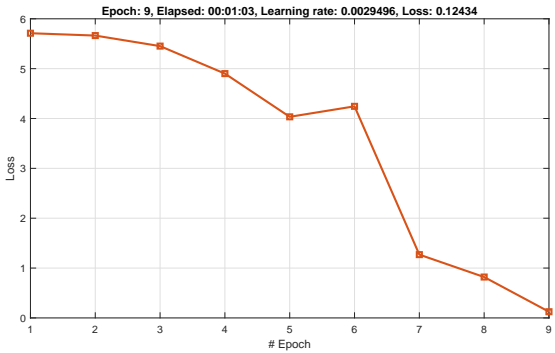}
    \caption{Training curve for example \#2.}
    \label{fig:trainingcurve}
\end{figure}

%
%
We consider the case where $t_0 = 0, \, T = 1$ and
\begin{equation}
\begin{aligned}
    &\Omega_0 = \{x\in\mathbb{R}^2 | x^T R x \le 1 \}, \, R = I \\
    &\Omega_t = \{x\in\mathbb{R}^2 | x^T \Gamma(t) x < 1 \},\, \Gamma(t) = 0.8k\mu(t)I \,.
\end{aligned}
\end{equation}
with $r_0=1$. From~\eqref{eq:laksh_bound} it can be seen that system~\eqref{eq:laksh_ex} is FTS with respect to $(t_0,T,\Omega_0,\Omega_t)$.
Some sample trajectories of system~\eqref{eq:laksh_ex} are shown in fig.~\ref{fig:ex2} together with the initial and trajectory domains.
With the parameters reported in tables~\ref{tab1}-\ref{tab2}, a solution is found after {9 training epochs (about $60$ s training)}.
The resulting $V(t,x)$ is shown at different time instants in fig.~\ref{fig:V3}, {while the resulting training curve (epochs vs. loss) for the considered parameters is shown in fig.~\ref{fig:trainingcurve}}.

{
In order to provide some insight on how the choice of the collocation points affects the algorithm performance, we performed a scan on the parameters $N_c, n_b, N_t, N_0$ of table~\ref{tab1} for this example. In particular, we varied the number of internal collocation points ($N_c = \{5, 10, 30, 50, 55, 70\} \times 10^{3}$) while keeping the other parameters fixed ($N_0 = 700, n_b = 200$ and a time step equal to $dt = 10^{-1}$, i.e. $N_t = 11$; see table~\ref{tab1}). Then, we set $N_c = 5\times10^{4}$ and varied the number of collocation points in the initial domain ($N_0 = 100$), the number of collocation points on the trajectory domain boundary for each time step ($n_b = \{50, 1000\}$) and the number of time instants~$N_t$ (by reducing the time step to $dt = 10^{-2} s$ while keeping $n_b = 200$). The results of this analysis are reported in fig.~\ref{fig:comparison}. The algorithm converged in all the considered cases, but for the ones represented by dashed lines (corresponding to cases with very few internal or boundary collocation points) a false positive was obtained, i.e. the resulting neural Lyapunov function violated the required conditions on some points of the test set. The test set for this analysis was generated by considering a different set of internal collocation points, obtained by restricting a $50 \times 50$ equally spaced grid to the initial and trajectory domains, while keeping the same boundary points. The time needed for the algorithm to reach convergence is also reported for completeness; however, notice that no effort has been made to optimize the code performance. Moreover, it is also worth to remark that a certain amount of randomness is implicitly involved in the procedure (in the initialization of the network parameters and in the choice of the collocation points) and that the used platform is not Real-Time. The case reported in fig.~\ref{fig:trainingcurve} corresponds to the red line in~\ref{fig:comparison}.

    \begin{figure}
        \centering
        \includegraphics[width=\linewidth]{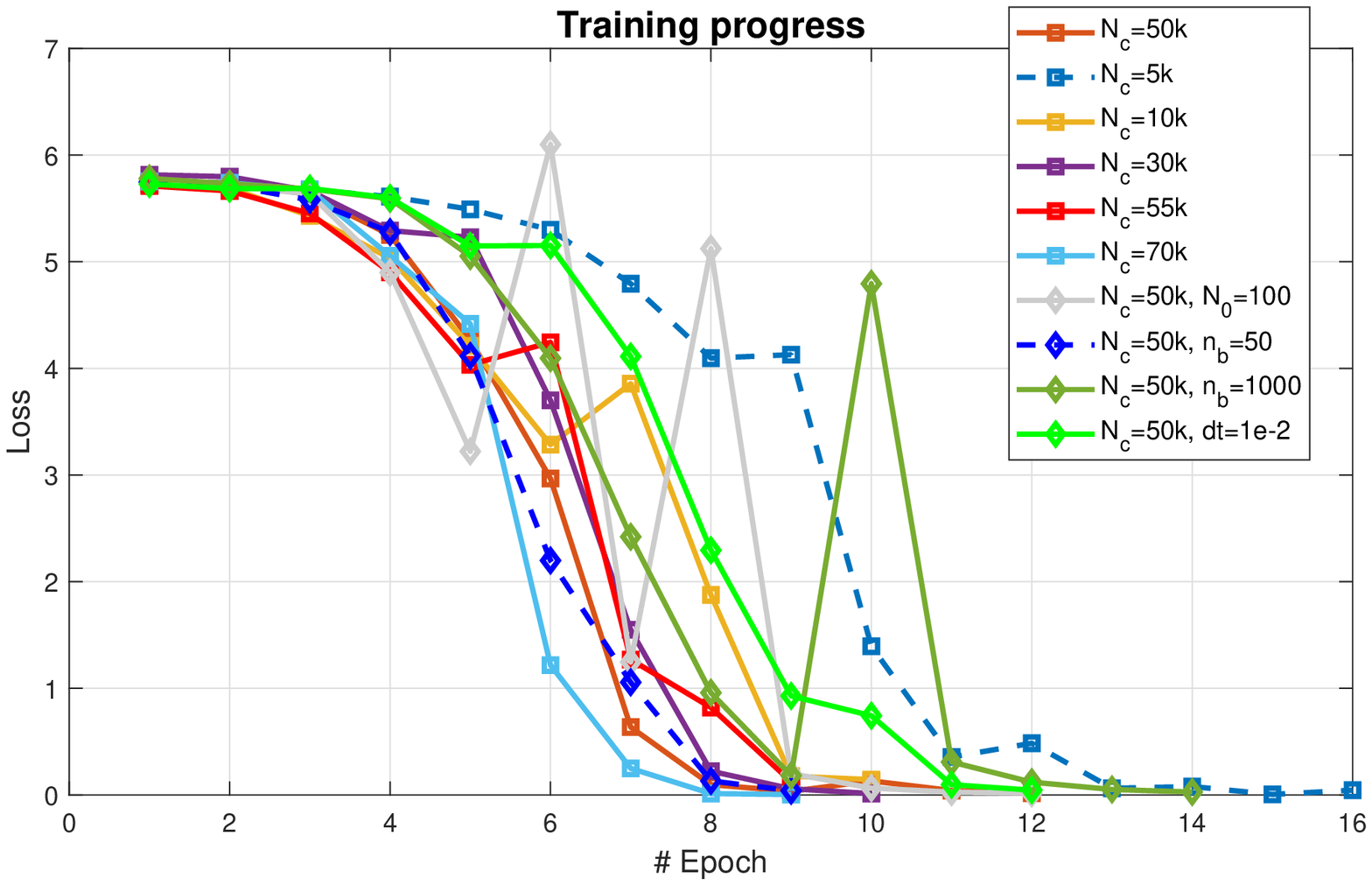}
        \includegraphics[width=\linewidth]{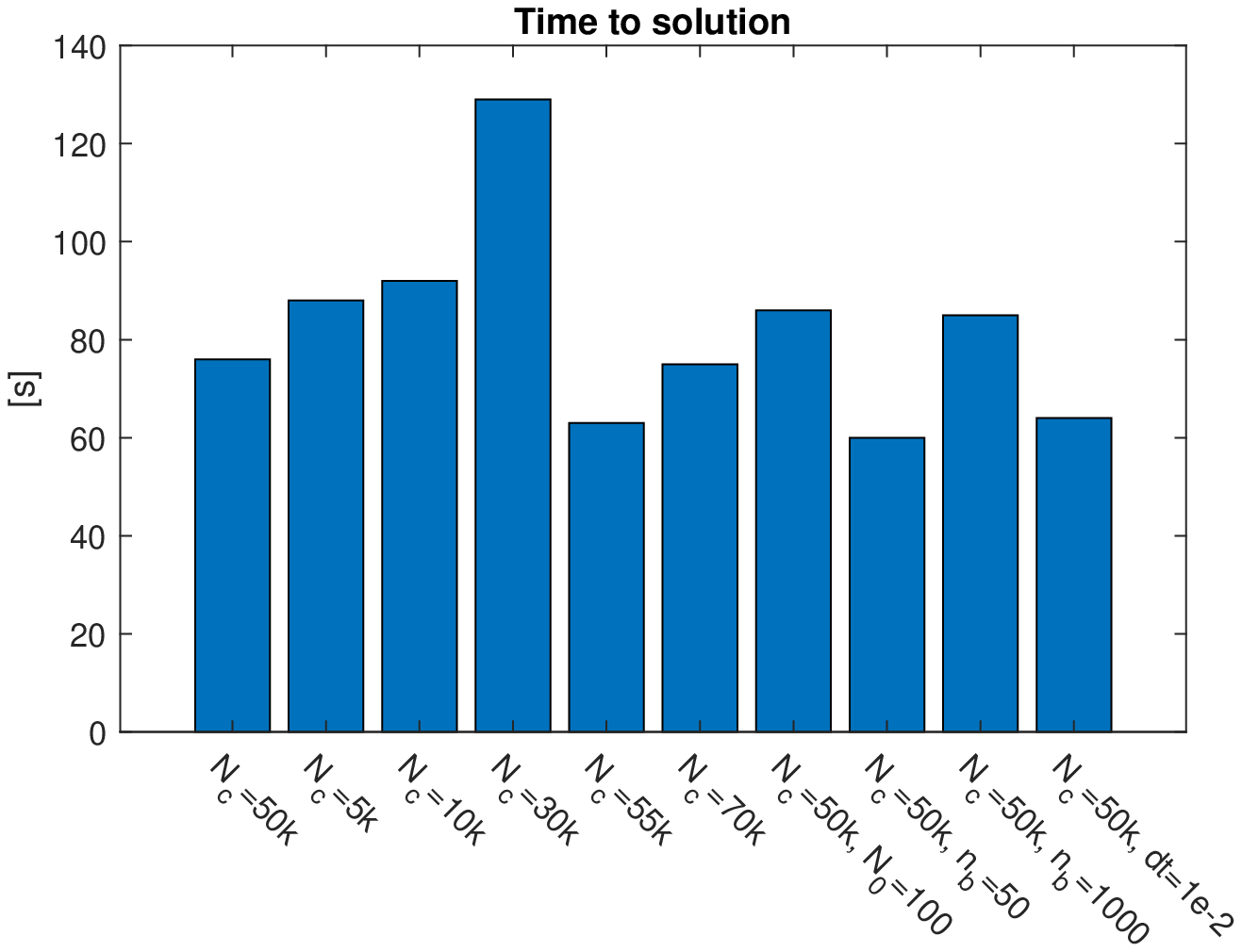}        
        \caption{Training curve for different choices of the parameters reported in table 1 for the case of example \#2 (top) and time to converge for each choice (bottom).}
        \label{fig:comparison}
    \end{figure}

    
}

\begin{figure}
    \centering
    \includegraphics[width=1\linewidth,trim=60 30 30 20]{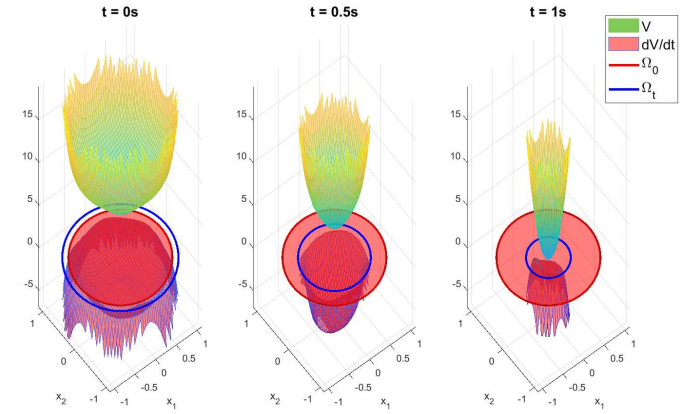}
    \caption{Resulting Lyapunov function $V(t,x)$ and its derivative $\dot{V}(t,x)$ for example \#2 at different time instants (a constant has been added so that the resulting $V(t,x)$ is above the $x_1-x_2$ plane for graphical reasons).}
    \label{fig:V2}
\end{figure}

\subsection*{Example 3} 
\begin{figure}[t]
    \centering
    \includegraphics[width=1\linewidth,trim=60 30 30 20]{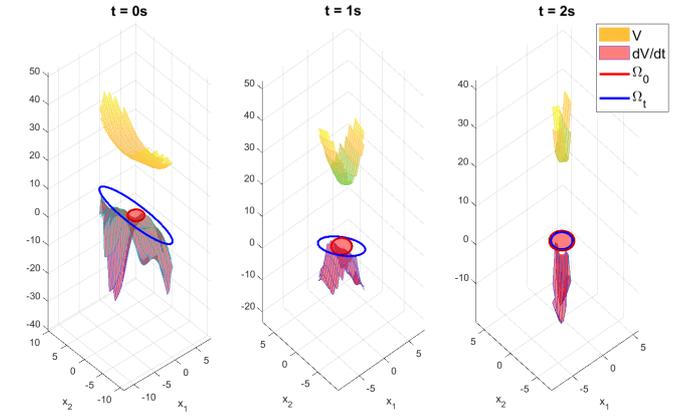}
    \caption{Resulting Lyapunov function $V(t,x)$ and its derivative $\dot{V}(t,x)$ for example \#3 at different time instants.}
    \label{fig:V3}
\end{figure}

\begin{figure}
    \centering
    \includegraphics[width=\linewidth]{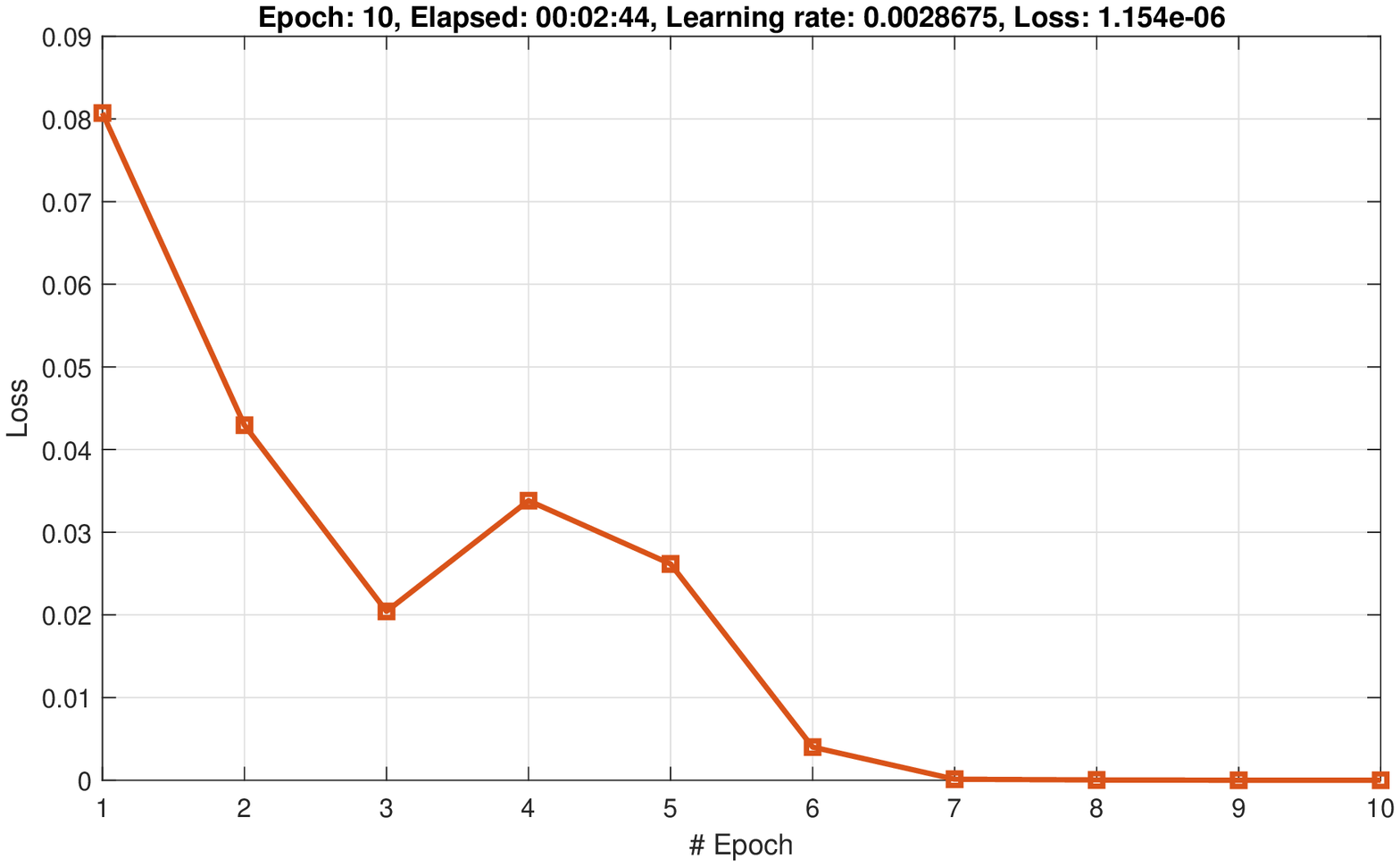}
    \caption{Training curve for example \#3.}
    \label{fig:trainingcurve2}
\end{figure}

\begin{figure}
     \centering
     \begin{subfigure}[b]{0.48\linewidth}
         \centering
         \includegraphics[width=0.9\textwidth, trim=20 0 20 0]{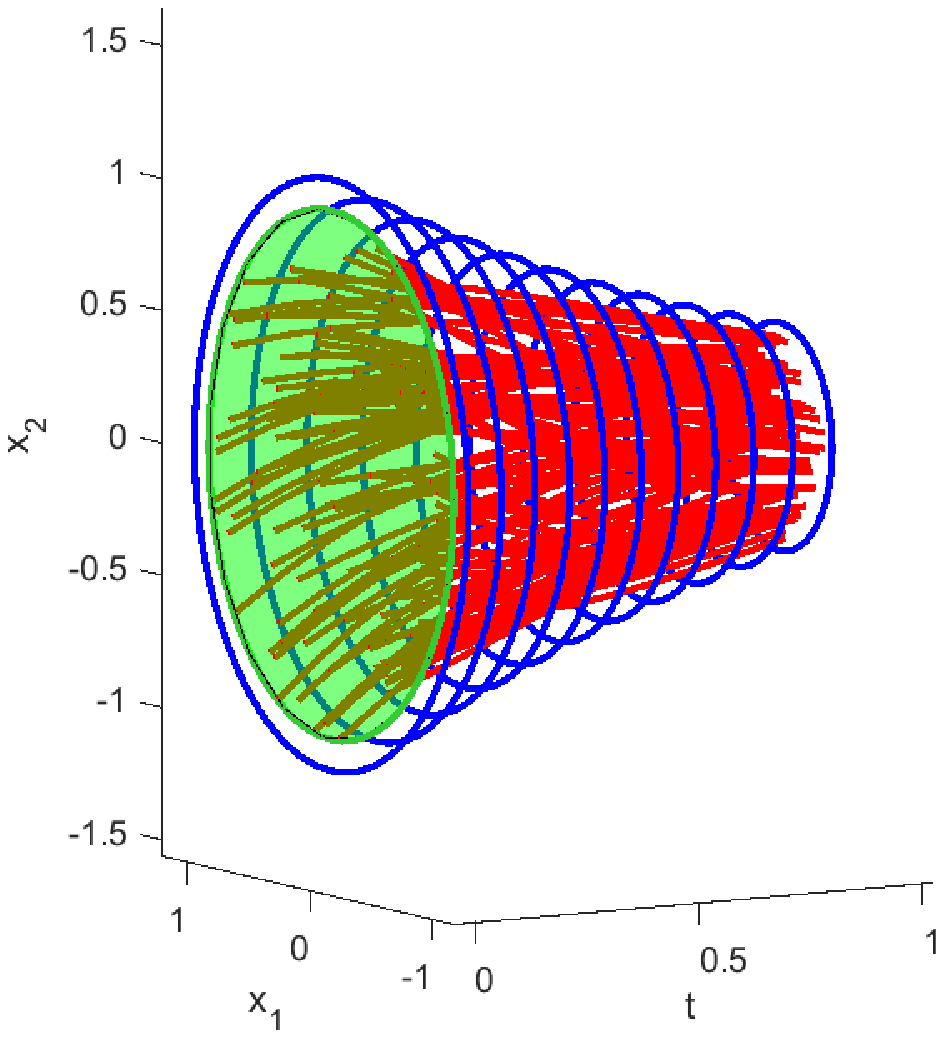}
         \caption{Trajectories of sys.~\eqref{eq:laksh_ex}.}
         \label{fig:ex2}
     \end{subfigure}
     \hfill
     \begin{subfigure}[b]{0.48\linewidth}
         \centering
         \includegraphics[width=0.9\textwidth, trim=20 0 20 0]{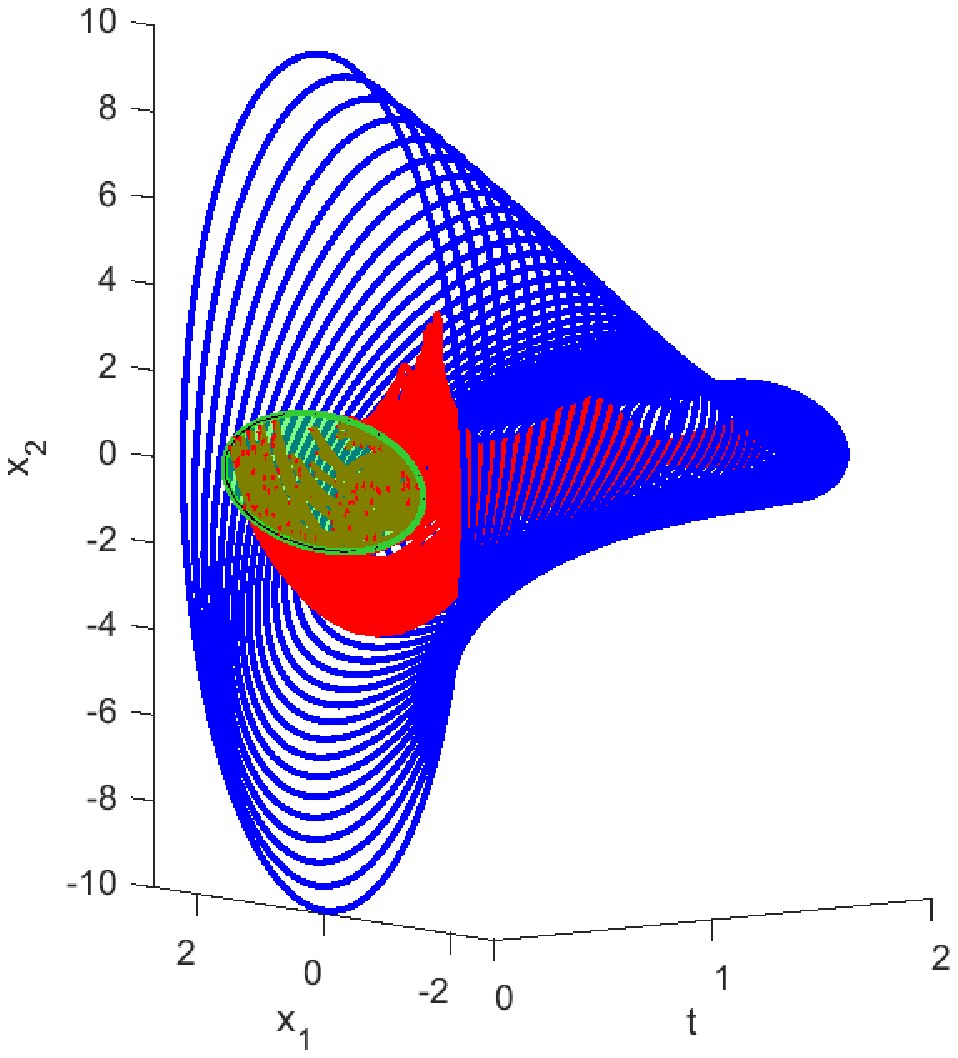}
         \caption{Trajectories of sys.~\eqref{eq:pendulum}.}
         \label{fig:ex3}
     \end{subfigure}
    \caption{Sample trajectories for the systems of examples \#2 and \#3 {(red)}, together with the considered initial (green) and trajectory (blue) domains.}
    \label{fig:trajectories}
\end{figure}

{Finally, consider the equations of a damped pendulum }
\begin{equation}\label{eq:pendulum}
\begin{aligned}
    &\dot{x}_1 = x_2\\
    &\dot{x}_2 = -\frac{g}{l}\sin(x_1) - \frac{b}{ml^2}x_2 \,,
\end{aligned}
\end{equation}
{which is a widely used benchmark problem in the field of systems and control theory. We set} $g=9.81$~m/s$^2$, $m=0.15$~kg, $l=0.15$~m, $b=0.1$~kg m$^2$/s. We consider the case where $t_0 = 0$, $T = 2$~s, and
\begin{equation}
\begin{aligned}
    &\Omega_0 = \{x\in\mathbb{R}^2 | x^T R x \le 1 \}, \, R = \frac{4}{\pi^2}I \\
    &\Omega_t = \{x\in\mathbb{R}^2 | x^T \Gamma(t) x < 1 \},\, \Gamma(t) = \Theta \tilde{\Gamma} \Theta^T \,,
\end{aligned}
\end{equation}
where
\begin{equation}
    \tilde{\Gamma} = \begin{bmatrix}
        \frac{2}{\pi^2} e^{0.5t} \quad 0 \\ 0 \quad \frac{0.1}{\pi^2}e^{2t}
    \end{bmatrix}
\end{equation}
and $\Theta$ is the rotation matrix given by
\begin{equation}
    \Theta = \begin{bmatrix}
        \cos(0.2\pi t) \quad \sin(0.2\pi t) \\ -\sin(0.2\pi t) \quad \cos(0.2\pi t) \,.
    \end{bmatrix}
\end{equation}
With the parameters reported in tables~\ref{tab1}-\ref{tab2}, a solution is found after 10 training epochs (approximately $2$~m $44$~s). 
Some sample trajectories for the system~\eqref{eq:pendulum} are shown in fig.~\ref{fig:ex3}, while the resulting $V(t,x)$ at different time instants is shown in fig.~\ref{fig:V3}. {The obtained training curve (epochs vs. loss) for the considered parameters is shown in fig.~\ref{fig:trainingcurve2}}.


\begin{table}[th]
\begin{center}
\footnotesize
\begin{tabular}{|c | c c c c c c c c|} 
 \hline
 Ex. & $\alpha_1$ & $\alpha_2$ & $\delta_1$ & $\delta_2$ & $N_c$ & $n_b$ & $N_t$ & $N_0$ \\ [0.5ex] 
 \hline
 \textbf{\#1} & 1 & 1 & 1 & 0.1 & 5000 & 500 & 11 & 200   \\ 
 \hline
 \textbf{\#2} & 5 & 1 & 1 & 3 & 55000 & 200 & 11 & 700  \\
 \hline
 \textbf{\#3} & 5 & 1 & 0.1 & 1 & 60000 & 100 & 15 & 700  \\
 [1ex] 
 \hline
\end{tabular}
\end{center}
 \caption{Training parameters for the proposed examples.} \label{tab1}
\end{table}

\begin{table}[th]
\begin{center}
\footnotesize
\begin{tabular}{|c | c c c c c|} 
 \hline
 Ex. &  hidden      & neurons & mini-   & learn & decay \\ [0.5ex] 
     &  layers      &         & batches & rate  & rate  \\ [0.5ex] 
 \hline
 \textbf{\#1} & 1 & 128 & 100 & 0.01 & $10^{-5}$  \\ 
 \hline
 \textbf{\#2} & 3 & 32 & 200 & 0.003 & $10^{-5}$ \\
 \hline
 \textbf{\#3} & 3 & 32 & 500 & 0.003 & $10^{-5}$ \\
 [1ex] 
 \hline
\end{tabular}
\end{center}
 \caption{Network structures for the proposed examples.} \label{tab2}
\end{table}

\section{Conclusions}\label{sec:conclusions}
{
In this article, we proposed some necessary and sufficient conditions for the FTS of nonlinear, non-autonomous dynamical systems in both continuous and discrete time. For the continuous-time case, the proposed conditions have been proven to be a generalization of standard LMI-based conditions that are customarily used for linear systems over ellipsoidal domains. 

A novel technique has then been proposed to assess the FTS property. The main idea is to exploit the capability of NN to serve as universal functions approximators in order to find a closed-form expression for the required Lyapunov-like function. The effectiveness of the approach has been illustrated through different numerical examples in the case of continuous-time systems. To the best of our knowledge, this is the first time that a numerical technique allows to practically verify this property for such a large class of systems. 

As discussed in the literature review of the introduction, the practical advantages of the FTS notion have been known since long time ago. However, after some pioneering works, the research community that had gathered around this topic remained almost silent until practical (D)LMI-based conditions to assess FTS emerged, at least in the case of linear systems, renewing the research interest in the field. Over the years, such conditions have been extended to uncertain, hybrid, stochastic, time-delay linear systems and to special classes of nonlinear systems. 
In this article, we proposed an alternative way to assess the FTS property for the case of general nonlinear systems. In this perspective, we expect that the method discussed here can be further refined and extended to more challenging classes of dynamical systems as happened for the linear case.

Of course, the approximate nature of the numerical method  employed to verify the conditions of thm.~\ref{thm1} leads to some limitations. First of all, since a discrete set of collocation points is considered, the proposed numerical test is in fact only necessary: the required conditions on the $V(t,x)$ function may be violated outside the training set of collocation points, but inside the domain of interest. However, this issue is quite standard when neural Lyapunov methods are considered. 
One possibility to mitigate this issue is to test the conditions on a set of points that is different (possibly also larger) with respect to the considered collocation points; this is the approach adopted in this paper. 
Furthermore, emerging research lines exist that focus on ways to assess the safety of similar methods in the framework of classic Lyapunov stability (see for example~\cite{abate2020formal}).
Another limitation of the proposed method is that there is no way of directly assessing that a system is Finite-Time \textit{unstable}. When this is the case, the training procedure will usually fail, as it cannot find a Lyapunov-like function that satisfies the required conditions. 
A possible future line of research could be that of extending classical Lyapunov instability theorems to the field of Finite-Time stability analysis.

Finally, the availability of a Lyapunov function in closed-form, whose derivatives with respect to time and state variables can be simply computed through auto-differentiation techniques, paves the way to novel Finite-Time Stabilization and Control techniques based on Lyapunov methods.
}




\bibliographystyle{model1-num-names}
\bibliography{references.bib}

\end{document}